\documentclass[a4paper,10.2pt]{article}

\usepackage{amsmath, amssymb}
\usepackage{geometry,amsmath,amssymb,enumerate,bbm,amsthm}
\usepackage[pdftex,bookmarks,colorlinks,breaklinks]{hyperref}  

\usepackage{memhfixc}  
\usepackage{pdfsync}  

\usepackage{amsfonts}
\usepackage{amsmath,amsthm,indentfirst}
\usepackage{amssymb}
\usepackage{amsthm}
\usepackage{amssymb}
\usepackage{amsmath}
\usepackage{vmargin}
\usepackage{dsfont}
\usepackage{mathrsfs}
\usepackage{fancyhdr}

\newcommand{\ep}{\epsilon}
\everymath{\displaystyle}

\newcommand{\RR}{\mathbb{R}}
\newcommand{\ZZ}{\mathbb{Z}}

\newtheorem{theo}{Theorem}
\newtheorem*{theo'}{Theorem 8 (Resnick)} 
\newtheorem*{theo''}{Theorem 9 (Abidi and Hmidi)}
\newtheorem*{lemmax}{Lemma 14 (maximum principle)}
\newtheorem{prop}[theo]{Proposition}
\newtheorem{lemm}[theo]{Lemma}

\newtheorem{rema}[theo]{Remark}
\newtheorem*{remak}{Concluding remarks}
\newtheorem{defi}[theo]{Definition}

\date{}
\author{Omar Lazar \thanks{  Universit\'e Paris-Est, LAMA (UMR 8050), UPEMLV, UPEC, CNRS, F-77454, Marne-la-Vall\'ee, France. \hspace{1cm} Email: omar.lazar@univ-mlv.fr}}
\title{Global existence for the critical dissipative surface quasi-geostrophic equation } 
\begin{document}
\maketitle
\bibliographystyle{plain}
\noindent{\bf Abstract:} In this article, we study the critical dissipative surface quasi-geostrophic equation (SQG) in $ \mathbb{R}^2$. Motivated by the study of the homogeneous statistical solutions of this equation, we show that for any large initial data $\theta_{0}$
liying in the space $\Lambda^{s} (H^{s}_{uloc}(\mathbb{R}^2)) \cap L^\infty(\mathbb{R}^2)$ the critical (SQG) has a global weak solution in time for $1/2<s<1$. Our proof is based on an energy inequality verified the  $(SQG)_{R,\ep}$ equation which is nothing but the equation (SQG) equation with truncated and regularized initial data. By classical compactness arguments, we show that we are able to pass to the limit ($R \rightarrow \infty$, $\ep \rightarrow 0$) in $(SQG)_{R,\ep}$ and that the limit solution has the desired regularity.

\section{Introduction}

We consider the initial value problem of the 2D dissipative surface
quasi-geostrophic equation :
\begin{equation}~\label{qg}
\ \\\left\{
\aligned
&\partial_{t}\theta(x,t)+u.\nabla\theta+ \Lambda^{\alpha}\theta = 0
\\ \nonumber
& u(\theta)= \mathcal{R}^\perp \theta \hspace{3,5cm}{:(SQG)_{\alpha}}
\\ \nonumber
& \theta(0,x)=\theta_{0}(x)
\endaligned
\right.
\end{equation}
where  \[\Lambda^{\alpha}\theta= (-\Delta)^{\alpha/2} \theta=C_{\alpha}P.V.\int_{\mathbb{R}^{n}}{\frac{\theta(x)-\theta(x-y)}{|y|^{2+\alpha}}dy},\]
$C_{\alpha}$ is a positive constant and $\theta\,:\,\mathbb R_{+}\times \mathbb R^2\to \mathbb R$ is a scalar function. Note that when  $0<\alpha<1$  and $\theta$ regular enough we can get rid of the principal value since the integral is then absolutely convergent.  Here $\alpha \in (0,2]$ is a fixed parameter and the velocity
$u=(u_1,u_2)$ is divergence free and determined by the Riesz
transforms of the potential temperature $\theta$ via the formula:
$$
u=(-{\mathcal R}_2\theta,{\mathcal R}_1\theta)=(-\partial_{x_2}
(-\Delta)^{-1/2}\theta,
\partial_{x_1}(-\Delta)^{-1/2}\theta).
$$
We can distinguish 3 cases depending on the value of $\alpha$. Namely, when $\alpha > 1$, $\alpha= 1$, $\alpha < 1$ which are called respectively the sub-critical, critical and super-critical cases. Actually, the more physically relevant case  is the critical one due to the term $(-\Delta)^{1/2}$ which models the Ekman pumping; as such it provides understanding of the quasi geostrophic flows. Moreover, some authors (see e.g. {\cite{CMT}}) have noticed that there is a similarity between the critical case and the 3D Navier-Stokes equation. In this paper, we will focus only on the critical case ($\alpha= 1$)  which we will denote by $(SQG)$.  The existence of weak solutions for this equation was proved by Resnick in \cite{Resnick} when the initial data lies in $L^{2}(\mathbb R^2)$. Many authors have studied the critical case, among them we can cite the work by Constantin, C\'ordoba and Wu {\cite{CCW}} in which they showed that there is a unique global solution when $\theta_{0} \in H^{1}$ and under a smallness assumption on $\Vert \theta_{0} \Vert_{L^\infty}$. In \cite{CV}, Caffarelli and Vasseur proved the global regularity of weak solutions with $L^{2}$ data. Their proof is based on De Giorgi techniques. Another important result is the one of Kiselev, Nazarov and Volberg {\cite{KNV}} wherein the authors used a non local maximum principle verified by the modulus of continuity at time $0$. They showed that all smooth periodic initial data gives rise to a unique global smooth solution. Abidi and Hmidi  \cite{Hm}  proved the global well-posedness  of the critical  dissipative quasi-geostrophic equation for large initial data  belonging to the critical Besov space $\dot B^0_{\infty,1}(\mathbb R^2)$. In \cite{Marchand}, Marchand showed the existence and regularity of global weak solutions to the quasi-geostrophic equations when 
the initial data belongs to $L^p$ $(p>4/3)$ or to $\dot H^{-1/2}$. \\

Whereas the sub critical case ($\alpha>1$) is globally well understood (see e.g \cite{CW}), the global regularity issue for large data in the super critical case ($\alpha<1$) is still open. Nevertheless, in the super critical case, several  results have been obtained, for instance, we can cite the work  by  Chae and Lee {\cite{CHAELEE}} where the authors showed the global regularity for small initial data in spaces $B_{2, 1}^{2 - 2 s}$. Global results for small initial data in Sobolev spaces $H^{m}$ with  $m<2$ can also be found in \cite{CC} and \cite{JU}.    \\

 In this paper, we show the existence of a global weak solution of $(SQG)$ for an initial data $\theta_{0}\in \Lambda^{s} ( H^{s}_{uloc}(\mathbb{R}^2)) \cap L^\infty(\mathbb{R}^2)$.  This space is a well adapted space when we plan to construct space-time homogeneous statistical solutions in the sense of Vishik and Fursikhov (\cite{VF1}, \cite{VF2}) for the $(SQG)$ equation. Roughly speaking, those statistical solutions have necessarily infinite energy, that is why the classical case  $\theta_{0}\in L^{2}$ is not allowed. Unlike the Navier-Stokes equation case, we cannot consider the case $\theta_{0}\in L^{2}_{uloc}$ as in the pioneer work by Lemari\'e-Rieusset \cite{PGLR}, \cite{PGLR2}, and more recently Basson \cite{Basson}.  This is due to the fact that Riesz transforms   are not well-defined in that space. Actually, to overcome this difficulty, we need more decay at infinity for the kernel of the Riesz transforms operator, so if we give more oscillations to $\theta_{0}$ we will be able to consider the case $\theta_{0} \in L^{2}_{uloc}$. As a matter of fact, if we put $\theta_{0} = \Lambda^{s}w_{0} \in \Lambda^{s} (H^{s}_{uloc}(\mathbb{R}^2))$  with $1/2<s<1$, then we can show that the $(SQG)$ equation makes sense in $\mathcal{D'}( \mathbb R^2)$. Since we plan to use the maximum principle we futher suppose that $\theta_{0} \in L^{\infty}$. \\

Our main result is the following theorem :

\begin{theo} 
Assume that $1/2< s<1$, and $\theta_0\in \Lambda^{s}(H^{s}_{uloc}) \cap L^{\infty},$  then the (SQG) equation has at least one global weak solution. Moreover for any $T<\infty,$  the solution belongs to the following spaces:
$$
\theta \in L^\infty ([0,T], L^{2}_{uloc}) \cap (L^{2}_{t}([0,T], \dot H^{1/2}))_{uloc} $$
and 
$$
w \in L^\infty([0,T], H^{s}_{uloc}) \cap (L^{2}_{t}([0,T], \dot H^{s+1/2}))_{uloc}.
$$
Futhermore, the following inequality holds for all $t\leq T$
$$
\Vert w(x,t) \Vert^{2}_{L^{\infty}_{t}{ H^{s}_{uloc}(\mathbb{R}^2)}} \leq c \ e^{CT},
$$

where $c$ and $C$ are two constants depending only on  $\displaystyle\Vert \theta_{0} \Vert_{L^{\infty}(\mathbb{R}^2)}$ and \ $\Vert w_{0} \Vert_{ H^{s}_{uloc}(\mathbb{R}^2)}$ .
\end{theo} 
\begin{rema}
 Beside giving an extension of Resnick's theorem (see theorem \ref{res}), this theorem would also allow us to show the existence of statistical solutions in the sense of Vishik and Fursikhov \cite{VF1}, \cite{VF2}. 
\end{rema}

This paper is organized as follows. In the next section we give a definition of the spaces and fix some notations. In the third section, we review some well known results we will use in the proof of our main theorem and we prove an energy inequality verified by the truncated equation associated with a regularized and truncated initial data. In the fourth section, we show that we can pass to the weak limit with respect to the parameters, this allows us to conclude the proof of our main result.

\section{Preliminaries and notations}

In this section, we recall some definitions and well-known results that we will use in the proof of our main theorem.\\

We start with the definition of the space $L^{p}_{uloc}(\mathbb R^2)$. This space belongs to the more general class of the so-called non homogenous Morrey-Campanato spaces defined by:
$$
M^{p}_{q}(\mathbb{R}^2)=\left\{ f \in L^{p}_{loc}(\mathbb{R}^{2}) \mid \sup_{x\in \mathbb{R}^2} \sup_{0<R<1} R^{\frac{2}{p}-\frac{2}{q}} \Vert f\Vert_{L^{q}(B(x,R))} < \infty \right\}.
$$
For $p=q$ we get the $L^{p}_{uloc}(\mathbb R^2)$ space.  In the sequel, we will use the  following definitions of the  $L^{p}_{uloc}(\mathbb R^2)$ spaces  :
\begin{defi}
Let us fix a positive test function $\phi_{0}$ such that $ \phi_{0} \in \mathcal{D} (\mathbb R^2) $ and  
\begin{equation}
\left \{
\aligned
&\phi_{0}(x)=1\ \  \mathrm{if} \ \vert x\vert \leq 2&
\\ \nonumber
& \phi_{0}(x)=0\ \ \mathrm{if} \  \vert x \vert \geq 3.&
\endaligned
\right.
\end{equation}
\end{defi}
We define the set of all translations of the function $\phi_{0}$ as $B_{\phi_{0}}=\{\phi_{0}(x-k), k\in \mathbb Z^2 \}.$ In the following, this set will be denoted  by $B$ or $B_{\phi_{0}}.$ 

\begin{rema}
Let $1\leq p \leq \infty$,  $f \in L^{p}_{uloc}(\mathbb R^2)$   if and only if  $f \in L^{p}_{loc}(\mathbb R^2)$ and the following norm is finite
 $$
\Vert f \Vert_{L^{p}_{uloc}(\mathbb R^2)}= \sup_{\phi \in B_{\phi_0}} \Vert \phi f \Vert_{L^{p}(\mathbb R^2)}.
$$
We will also use the following useful equivalent norms:
$$\Vert f \Vert^{p}_{L^{p}_{uloc}(\mathbb R^2)} \approx \sup_{k \in \mathbb{Z}^2 } \int_{k+[0,1]^2} \vert f(x) \vert^p \ dx \approx \sup_{k \in \mathbb{Z}^2 } \Vert \phi (x-k) f \Vert_{L^{p}(\mathbb R^2)}  $$
\end{rema}

Let us also recall the following useful lemma:
\begin{lemm}
Let $f\in L^{1}(\mathbb{R}^2)$ and $g\in L^{p}_{uloc}(\mathbb R^2)$ then the convolution product is well defined $(a.e.)$ and we have $\Vert f * g \Vert_{L^{p}_{uloc}(\mathbb R^2)} \leq \Vert f \Vert_{L^{1}} \Vert g \Vert_{L^{p}_{uloc}}.$
\end{lemm}
\noindent{\bf Proof.}
As mentioned in \cite{Basson}, if we put $K=x_{0}+[0,1]^2$, it suffices to write:
$$
\left\Vert \int \vert f(y) g(x-y) \vert \ dy \right\Vert_{L^{p}(K)} \leq \int \vert f(y) \vert \Vert g(x-y) \mathds{1}_{x \in K}\Vert_{L^{p}_{x}} \ dy \leq \Vert f \Vert_{L^{1}} \Vert g \Vert_{L^{p}_{uloc}}
$$

\qed

Throughout this paper, we will make  use of a test function $\psi$ constructed as follows. We introduce a positive test function $\psi_{0}$  supported in $[-4,4]^2$ such that $\psi_0 \geq 0$ and 
\begin{equation}
\left \{
\aligned
&\phi_{0}(x)=1\ \  \mathrm{if} \ \vert x\vert \leq 3&
\\ \nonumber
& \phi_{0}(x)=0\ \ \mathrm{if} \  \vert x \vert \geq 4.&
\endaligned
\right.
\end{equation}

Now, we construct the space of all translations of $\psi_{0}$ that is:
$$
B_{\psi_{0}}=\{ \psi_{0}(x-x_{0}), \ x_{0} \in \mathbb{Z}^2 \}.
$$

Then, as we did before, we define $\psi(x)=\psi_{0}(x-k)$ where $k \in \mathbb{Z}^2$. Note that by construction, the function $\psi$ is equal to 1 in a neighborhood of the support of $\phi$. \\

We will also denote $A \lesssim B$ if $A$ is less than $B$ up to a positive multiplicative constant which can be different from line to line, and the positive constant $C$ which appears in some estimations can be different as well. Keep in mind that those constants may depend on certain controlled norms of the initial data. \\

We recall the definition of the Sobolev space $H^{s}(\mathbb{R}^{2})$ where $s\in \mathbb{R}$ :

\begin{defi}
A distribution $f \in \mathcal{S}'(\mathbb{R}^{2})$ is in $H^{s}(\mathbb{R}^{2})$ if $\hat{f}$ is locally integrable on $\mathbb{R}^{2}$ and
\begin{equation*}
\left\Vert f\right\Vert _{H^{s}}\equiv \left( \int_{\mathbb{R}^{2}}\left(
1+\left\vert \xi \right\vert ^{2}\right) ^{s}\left\vert \hat{f}(\xi
)\right\vert ^{2}d\xi \right) ^{1/2} \ \ <\infty .
\end{equation*}
\end{defi}

The Sobolev spaces can be defined through the fractional Laplacian, and more precisely we have the following lemma.

\begin{lemm}
If $0<s<1$ then the norm :
\begin{eqnarray*}
N(f) &=& C_s \left( \int \vert f(x) \vert^{2} \ dx + \int \int \frac{\vert f(x)-f(y) \vert^2}{\vert x-y \vert^{n+2s}} \ dx \ dy \right)^{1/2}
\end{eqnarray*}
is finite on ${H}^{s}$ and this norm is equivalent to the usual ${H}^{s}$ norm. 
\end{lemm} 
\noindent{\bf Proof.} See e.g. Prop 1.37 p 28, \cite{BCD} \\


Now, let us recall the definition of the $H^{s}_{uloc}(\mathbb R^2)$ space, where $s\in \mathbb R $ :
\begin{defi}
Let $\phi_{0}$ be a positive test function such that $ \phi_{0}  \in \mathcal{D} (\mathbb R^2) $ such that $\phi(x)= \ 1 \ \mathrm{if} \ x \in [-1,1]^2$  and let $B=\{\phi_{0}(x-k), k \in \mathbb Z^2 \} $. We say that $f \in H^{s}_{uloc}(\mathbb R^2)$   if and only if  $f \in H^{s}_{loc}(\mathbb R^2)$ and if the following norm is finite:
\begin{center}
$\displaystyle \Vert f \Vert^{2}_{{H^{s}_{uloc}} (\mathbb R^2)}  = \sup_{\phi\in B  } \Vert \phi f \Vert_{H^{s}}$
\end{center}

\end{defi}
\begin{rema}
We can easily see that those defintions do not depend on the choice of the test function since we take the supremum over all translations of $\phi$. Let us precisely state that $f\in \Lambda^{s}\left.( \dot H^{s}_{uloc}(\mathbb R^2)\right.)$ if there exists $g\in \dot  H^{s}_{uloc}(\mathbb R^2)$ s.t $f=\Lambda^{s} g$. Let us also define some space-time norms  that we will use throughout this paper, namely $(L^{\infty}_{T} \dot H^{s})_{uloc}$, $(L^{2}_{T}\dot H^{s})_{uloc}$, where $s \in [0,1]$.
\begin{eqnarray*}
\Vert w \Vert^{2}_{ (L^{2}_{T} \dot H^{s})_{uloc}}&=& \sup_{\phi \in B} \int_{0}^{T} \int \phi \vert \Lambda^{s} w(x,s) \vert^2 \ dx \ ds < \infty, \\
\Vert w \Vert_{(L^\infty_{T} \dot H^{s})_{uloc}} &=& \sup_{t\in[0,T]} \sup_{\phi \in B} \int \phi \vert \Lambda^{s} w(x,t) \vert^2 \ dx  < \infty. 
\end{eqnarray*}
\end{rema}


In the sequel, we shall use the following notations: 
\begin{center}
\hspace{-1cm}$A_{\phi}(w)=\displaystyle \int \left(\frac{\vert w \vert^2}{2}\phi+\frac{\vert \Lambda^{s}w\vert^2}{2}\phi\right) dx $ \ \   and \ \  $A(w)=\displaystyle\int \left(\frac{\vert w \vert^2}{2}+\frac{\vert \Lambda^{s}w\vert^2}{2}\right) dx$.
\end{center}

The next lemma will be useful for the proof of our result :
\begin{prop}
Let $0<s<1$ ; the following norms are equivalent on $ H^{s}_{uloc}(\mathbb{R}^2)$   :
\begin{center}
\begin{itemize}
\item $ \ \Vert w \Vert_{a} = \sup_{\phi\in B}\{ A_{\phi}(w)< \infty \}$ 
\item $ \ \Vert w \Vert_{b} = \sup_{\phi\in  B}  \{A(w\phi)< \infty \}$. 
\end{itemize}
\end{center}
\end{prop}
\noindent{\bf Proof.}

 Since $0\leq \phi \leq 1$, we have $\frac{\vert w\phi \vert^2}{2} \leq    \frac{\vert w \vert^2}{2}\phi$ \ and then
 $$
 \sup_{\phi \in B} \int \frac{\vert w\phi \vert^2}{2} \ dx  \leq    \sup_{\phi \in B} \int \frac{\vert w \vert^2}{2}\phi \ dx. 
 $$
 
 Now, let us prove the reverse inequality. Suppose that we control $ \sup_{\phi \in B} \int \frac{\vert w \vert^2}{2}\phi \ dx$. 
 We put $\eta^{2}(x)=\phi(x)$, then
 $$
  \sup_{\phi \in B} \int \frac{\vert w \vert^2}{2}\phi \ dx  \leq \sup_{\eta^{2} \in B_{\eta^{2}}} \int \frac{\vert w\eta \vert^2}{2} \ dx.
 $$
 These two previous inequalities show the following equivalence of norm on $L^{2}_{uloc}$ 
 
 $$\Vert w \Vert^{2}_{L^{2}_{uloc}} \sim  \sup_{\phi \in B} \int \frac{\vert w \vert^2}{2}\phi  \ dx \sim  \sup_{\phi \in B} \int \frac{\vert w\phi \vert^2}{2} \ dx. $$
 
  Now, we just have to prove that if $\Lambda^{s}(\phi w) \in L^2$ then $\phi \Lambda^{s}w \in L^2$ and the converse.  Throughout this paper, we will frequently make  use of the commutator between an operator $T$ and a function $\phi$  defined by the formula
  $$
 [T,\phi] f= T(\phi f)-f T(\phi).
  $$
  
 Then, for all $\phi \in B_{\phi_{0}} $ and $\psi \in B_{\psi_{0}} $ we can write
$$
\psi [\Lambda^{s}, \phi] w+(1-\psi)[\Lambda^{s}, \phi] w=\Lambda^{s}(\phi w)- \phi \Lambda^{s} w.
$$
  
  The following lemma allows us to finish the proof of the proposition.
  
  \begin{lemm}
  The operator $T_{\phi} w$ defined by 
\begin{eqnarray*}
   T_{\phi} : &L^{2}_{uloc}& \longrightarrow L^{2}_{uloc}  \\
  &w& \longmapsto [\Lambda^{s}, \phi] w
\end{eqnarray*} 
is continuous.
  \end{lemm}
  \noindent{\bf Proof.} We write
  \begin{eqnarray*}
\left\vert [\Lambda^{s}, \phi] w \right\vert  &\leq& C_s\int \frac{\vert \phi(x)-\phi(y)\vert}{\vert x-y \vert^{2+s}} {\vert w (y)\vert} \ dy \\
&\leq& C_s \int \frac{\min\left(\Vert \phi \Vert_{\infty}, \Vert \nabla\phi \Vert_{L^\infty} \vert x-y \vert \right) }{\vert x-y \vert^{2+s}}{\vert w (y)\vert} \ dy \\
&\leq& C_s \int \frac{\min\left(1, \vert x-y \vert \right)}{\vert x-y \vert^{2+s}}{\vert w (y)\vert} \ dy \\
&=& (A* \vert w \vert)(x),
\end{eqnarray*}
where $A(x)= C_s \  \frac{\min\left(1, \vert x\vert \right)}{\vert x \vert^{2+s}} $. Then, we see that if $ \vert x \vert < 1 $  then $A(x)=\frac{C_s}{\vert x \vert^{1+s}} \in L^1$ and since $\vert w \vert \in L^{2}_{uloc}$ by convolution (see Lemma 2) we get that  $A* \vert w \vert \in L^{2}_{uloc}$. \\ If \ $\vert x \vert > 1 $ then $A(x)=\frac{C_s}{\vert x \vert^{2+s}} \in L^1$ and the same conclusion holds and this finishes the proof of the lemma.  

Thus, $\psi [\Lambda^{s}, \phi] w \in L^2$. For the term $(1-\psi)[\Lambda^{s}, \phi] w$ we see that:
\begin{eqnarray*}
(1-\psi)[\Lambda^{s}, \phi] w &=& (1-\psi) \Lambda(\phi w)-(1-\psi) \phi \Lambda w \\
&=& (1-\psi)\left(C_{s}\frac{\mathds{1}_{x\geq1}}{\vert x \vert^{2+s}}* w\phi \right).
\end{eqnarray*}
The kernel $\frac{C_{s} \mathds{1}_{x\geq1} }{\vert x \vert^{2+s}}$ is in $L^1$, and since $w\phi \in L^2$,  by convolution, we get that $(1-\psi)[\Lambda^{s}, \phi] w \in L^2$ and this ends the proof of the proposition.\\

\qed


\section{The energy inequality}

\subsection{Well known existence results in the critical case.}
In this section, we recall some existence and 
uniqueness results that we will use in the proof of our main theorem. Let us recall that the existence of global weak solutions of $(SQG)$ with initial data in $L^{2}(\mathbb{R}^2)$ has been obtained by Resnick in \cite{Resnick}. More precisely, he proved the following theorem:
\newpage
 \begin{theo'} \label{res}
 Let $\theta_{0} \in L^{2}(\mathbb{R}^2)$. Then, for any $T > 0$, there exists at least one weak 
solution to the critical (SQG) equation in the following sense:
$$
\partial_{t}\int \theta \phi \ dx-\int \theta(u.\nabla \phi) \ dx+\int (\Lambda^{1/2}\phi)(\Lambda^{1/2}\theta) \ dx=0, \ \ \forall \phi \in \mathcal{D} $$
Moreover, 
$$
 \theta \in L^\infty([0,T],L^{2}(\mathbb{R}^2)) \cap L^{2}([0,T], \dot H^{1/2}(\mathbb{R}^2)).
$$
 \end{theo'}
 
From the result of Caffarelli and Vasseur \cite{CV}, those weak solutions are known to be smooth. Unfortunately, the smoothness is not sufficient in our proof since in our future computations we need a little bit of integrability (at least $ H^{1}$) therefore we need a result of global and regular enough solutions  without condition of smallness on the initial data. As we recalled in the introduction, Abidi  and Hmidi in \cite{Hm} proved a theorem of global solution for large data in the critical Besov space $\theta_0 \in \dot B^0_{\infty,1}$. Let us recall their main theorem; to do so, we need to recall the definition of the homogenous Besov spaces based on the classical Littlewood-Paley decomposition.

Let
$\phi$ be a smooth function  supported in the ring $\mathcal{C}:=\{ \xi\in\RR^2,\frac{3}{4}\leq|\xi|\leq\frac{8}{3}\}$ and such that 
$$
\sum_{q\in\ZZ}\phi(2^{-q}\xi)=1 \quad\hbox{for}\quad \xi\neq 0.
$$
Now,  for  $u\in{\mathcal S}'$ we set
  $$
\forall q\in\ZZ,\quad \Delta_qu=\phi(2^{-q}\textnormal{D})u\hspace{1cm}\mbox{and}\hspace{1cm}
S_qu=\sum_{j\leq q-1}\Delta_{j}u. 
 $$ 
Recall that we have the following formal Littlewood-Paley decompostion of $u$
$$
u=\sum_{q\in\ZZ}\Delta_q \,u,\quad\forall\,u\in {\mathcal {S}}'(\RR^2)/{\mathcal{P}}[\RR^2],
$$
where ${\mathcal{P}}[\RR^2]$ is the set of polynomials. \\

 Now we are able to give the  definition of the homogenous Besov spaces. Let $(p,m)\in[1,+\infty]^2,$ $s\in\RR$ and $u\in{\mathcal S}',$ $u$ is said to be in $ \mathcal { \dot B}^{s}_{p,m}$ if
$$
\|u\|_{\mathcal { \dot B}^{s}_{p,m}}:=\Big(2^{qs}\|\Delta_q u\|_{L^{p}}\Big)_{\ell ^{m}} <\infty.
$$

In \cite{Hm}, Abidi and Hmidi proved the following regularity theorem:

 \begin{theo''}
Let $\theta_0 \in \dot B^0_{\infty,1},$ then there exists a unique global solution $\theta$ to the dissipative critical SQG equation such that
$$
\theta\in {\mathcal{C}}(\RR_+;\,\dot B^0_{\infty,1})
\cap L^1_{\textnormal loc}(\RR_+;\,\dot B^1_{\infty,1}).
$$
\end{theo''}

\begin{rema} The space $\dot B^0_{\infty,1}$ is defined as the completion in the Schwartz class for the norm $L^\infty$  with respect to the norm $\|u\|_{\mathcal { \dot B}^{0}_{\infty,1}}=\sum_{q\in \mathbb{Z}}\|\Delta_q u\|_{L^{\infty}} $.  
\end{rema}

We recall the  well-known result due to A. C\'ordoba and D. C\'ordoba {\cite{CC}}.  
\begin{lemm}
Let  $0 \leq \alpha \leq 2 $, the following pointwise inequality holds for all convex functions $g$, 
\begin{center}
$\Lambda^{\alpha}(g(\theta))\leq g'(\theta) \Lambda^{\alpha}(\theta)$
\end{center}
\end{lemm}
\noindent{\bf Proof.} See e.g.  {\cite{CC}}.\\

We shall use the maximum principle for the (SQG)  equation due to Resnick \cite{Resnick} and also A.C\'ordoba and D.C\'ordoba {\cite{CC}}, that is:
\begin{lemmax}  \
Let $\theta$ be a smooth function  satisfying the $(SQG)_\alpha$ equation \\ $\partial_{t}\theta+u.\nabla\theta+ k\Lambda^{\alpha}\theta = 0$, where $0\leq \alpha \leq 2$ and $k\geq 0$, then we have for all $p \in [2, +\infty)$ and for all $t \geq 0$
$$ 
\Vert \theta\Vert^{p}_{L^p}+2k\int^{t}_{0}\int \vert \Lambda^{\alpha/2}(\vert \theta\vert^{p/2}) \vert^2 \ dx \ ds \leq \Vert \theta_{0}\Vert^{p}_{L^p}
$$
and
$$
\Vert \theta(t)\Vert_{L^\infty}\leq \Vert \theta_{0}\Vert_{L^\infty}
$$
\end{lemmax} 
\noindent{\bf Proof.} See e.g.  {\cite{CC}}.

\subsection{The critical SQG equation with $\theta=\Lambda^{s}w$ }
Let us recall that we want to give some oscillations to $\theta$, so we put  $\theta =\Lambda^{s}w=(-\Delta)^{s/2}w$ and we suppose that $w\in   H^{s}_{uloc}$ and that $\theta\in L^{\infty}$. The (SQG) equation becomes:
\begin{equation*}
\Lambda^{s}(\partial_{t}w)=\nabla\cdot \left(\Lambda^{s}w (\mathcal{R}^{\perp}{\Lambda^{s} w})\right) -\Lambda^{s+1} w, 
\end{equation*}

hence, we get
\begin{equation}
\partial_{t}w=\left(\Lambda^{-s}\nabla \right)\cdot \left(\Lambda^{s}w (\mathcal{R}^{\perp}{\Lambda^{s} w})\right) -\Lambda{w} 
\end{equation} 

The equation that we will study through this paper is the next one:

\begin{equation}
(\stackrel{\sim}{SQG})_{} : \ \\\left\{
\aligned
&\partial_{t}w=\left(\Lambda^{-s}\nabla \right)\cdot \left(\Lambda^{s}w (\mathcal{R}^{\perp}{\Lambda^{s} w})\right) -\Lambda{w}
\\ \nonumber 
&   \nabla.\mathcal{R}^{\perp}{\Lambda^{s} w}= 0.
\endaligned
\right.
\end{equation}

Note that the initial value problem is now endowed with the conditions:
\begin{eqnarray*}
&&\theta_{0}(x)=\Lambda^{s} w_0 \in \Lambda^{s}(  H^{s}_{uloc}) \cap L^\infty.
\end{eqnarray*}


Before we go any further, we have to check that the $(\stackrel{\sim}{SQG})$ equation makes  sense  when  $w \in L^{\infty}([0,T], \dot H^{s}_{uloc}(\mathbb R^2))$ and $\Lambda^{s} w \in L^{\infty}([0,T],L^{\infty}(\mathbb R^2))$ where $T>0$. The equation $(\stackrel{\sim}{SQG})$ is:  

$$
\partial_{t}w=\Lambda^{-s}\nabla \cdot \Lambda^{s}w (\mathcal{R}^{\perp}{\Lambda^{s} w}) -\Lambda{w}
$$
If $w \in L^{\infty}([0,T], H^{s}_{uloc}(\mathbb R^2))$ then the computation $\Lambda w$ has obviously a sense in $\mathcal{D}'$ and $\Lambda^{s} w \in L^{\infty}([0,T], L^{2}_{uloc}(\mathbb R^2))$. We now deal with the non linear term which is:
\begin{center}
$\Lambda^{-s}\nabla \cdot \Lambda^{s}w (\mathcal{R}^{\perp}{\Lambda^{s} w})$
\end{center}
We begin with the study of the singular integral operator  $\mathcal{R}^{\perp}{\Lambda^{s}}$. Recalling that  $\mathcal{R}^{\perp} \theta = \nabla^{\perp} \Lambda^{-1} \theta$, we infer $\mathcal{R}^{\perp} \theta  =\nabla^{\perp} \Lambda^{s-1} w=(-\partial_{2}  \Lambda^{s-1} w, \partial_{1}  \Lambda^{s-1} w)$, since the convolution kernel of the operator $\Lambda^{s-1} \theta$  is defined by   $K(x)=\frac{C_s}{\vert x \vert^{1+s}}$; therefore, we see that the convolution kernel $\tilde K$ of $u_R$ is given by:
$$\tilde K(x)=\partial_{j} \frac{C_s}{\vert x \vert^{1+s}}=\frac{-C_s(1+s)}{\vert x \vert^{2+s}}  \frac{x_j}{\vert x \vert}.$$ 
Now, let $\alpha$ be a smooth cut-off function which is equal to one for $\vert x \vert \leq 1$ and zero for $\vert x \vert \geq 2$. Then we split,
$$
u_R=  \alpha \tilde K * w + (1-\alpha) \tilde K * w.
$$
Since $ \tilde K \alpha  \in \mathcal{E}'$ and $w \in L^{2}_{uloc} \subset   \mathcal{D}'$ the convolution makes sense. For the second part, we notice that $ (1-\alpha)  \tilde K  * \theta \in L^{1} * L^{2}_{uloc} \subset   \mathcal{D}'$. \\




  The idea of the proof of our main result is the following. We introduce a truncated initial data:
 $$\theta_{0,R}=\Lambda^{s}(w_{0}\chi_{R}),$$
  where $\chi_{R}$ is a positive smooth function constructed as follows. Let $\chi \in \mathcal{D}(\mathbb R^{2})$ be a positive smooth function s.t $\chi(x)=1$ if $\vert x \vert \leq 1$, and 0 if $\vert x \vert \geq 2$. For $R>0$, we introduce the function  $\chi_R (x)\equiv \chi(x/R)$.  By construction, the function $\chi_R $ is a positive smooth function s.t for all $R>0$,
 \begin{equation}
\\\left \{
\aligned
& \chi_{R}(x)=\ 1\  \mathrm{if}  \ \vert x \vert \leq R,&
\\ \nonumber
&  \chi_{R}(x)= \ 0 \ \mathrm{if} \  \vert x \vert \geq 2R.&
\endaligned
\right.
\end{equation}

We also need to introduce the following standard mollifier $\rho \in \mathcal{D}(\mathbb{R}^2)$, s.t $supp \rho \subset [-1, 1]^2$ and $\int_{\mathbb{R}^2} \rho = 1$, then we define $\rho_{\ep}(x)=\frac{1}{\ep^2}\rho (\frac{x}{\ep})$. \\
  
  We will study the following truncated equation associated with a truncated and regularized initial data. 
  
  \begin{equation}
\ \\\left\{
\aligned
&\partial_{t}w_{{R}}=\left(\Lambda^{-s}\nabla \right)\cdot \left(\Lambda^{s}w_{{R}}\  \mathcal{R}^{\perp}{\Lambda^{s} w_{{R}}}\right) -\Lambda{w_{{R}}}
\\ \nonumber
& \nabla.\ \mathcal{R}^{\perp}{\Lambda^{s} w_{{R}}}=0 \hspace{4,9cm} :({SQG})_{R,\ep} 
\\ \nonumber
& \theta_{0,R,\ep}=\Lambda^{s}(w_{0}\chi_{R}) * \rho_{\ep}
\endaligned
\right.
\end{equation}
  
Since the initial data $\theta_{0}\in \Lambda^{s} (H^{s}_{uloc}(\mathbb{R}^2)) \cap L^{\infty}(\mathbb{R}^2)$ then the truncated initial data $\theta_{0,R} \in L^{2 }(\mathbb{R}^2) \cap L^{\infty}(\mathbb{R}^2)$ and  
$$
\displaystyle\theta_{0,R, \ep} \in \displaystyle \bigcap_{k\geq 0} H^{k} \subset H^{2} \subset  \dot B^0_{\infty,1}
$$ 
Therefore, we can use the  result of Abidi and Hmidi \cite{Hm} and claim that for the truncated and regularized initial data, there exists at least one global weak solution $w_{R}$  to the $(SQG)_{R,\ep}$ equation. Moreover, the solution belongs to the following space:
 $$
\theta_{R} \in{\mathcal{C}}(\RR_+;\,\dot B^0_{\infty,1})
\cap L^1_{\textnormal loc}(\RR_+;\,\dot B^1_{\infty,1}).
  $$
  
Then, we will obtain an energy inequality involving the truncated solution and using some estimates and compactness arguments, we will show that we can pass to the weak limit of the $(SQG)_{R,\ep}$ equation. The solution of the $(SQG)$ equation will be the weak limit of $w_{R,\ep}$. \\

The following lemma will allow us to apply the maximum principle to $w_0$.

\begin{lemm}
If $w \in  H^{s}_{uloc}(\mathbb{R}^2)$ and $\theta=\Lambda^{s}w \in L^\infty(\mathbb{R}^2)$ then $w \in L^{\infty}(\mathbb{R}^2)$
\end{lemm} 

\noindent{\bf Proof.} We decompose the operator $\Lambda^{s}$ in low frequencies and high frequencies. We denote $S_{0}$ the low frequencies and we write:
$$
w=\Lambda^{-s}(Id-S_{0})\Lambda^{s}w+S_{0}w.
$$
Since the operator $\Lambda^{-s}(Id-S_{0})$ is continuous from $L^\infty$ to $L^\infty$ thus $\Lambda^{-s}(Id-S_{0}w)\Lambda^{s}w \in L^\infty$. Moreover, since  $ H^{s}_{uloc}(\mathbb{R}^2)$ is a shift invariant space, one can easily see that $S_{0}w \in L^\infty$.\\

\qed
\subsection{Bounds for the truncated and regularized initial data}
Now, we have to check that the truncated and regularized initial data lies in $L^\infty$ uniformly with respect to the parameters $R$ and $\ep$.

\begin{lemm}
If  $\theta_{0} \in L^\infty$ and $w_{0} \in  H^{s}_{uloc}(\mathbb{R}^2)$ then $\sup_{R>1,\ep>0} \Vert \theta_{0,R,\ep}  \Vert_{L^{\infty}(\mathbb{R}^2)} < \infty$
\end{lemm}
\noindent{\bf Proof.} Using a commutator, one can write 
$$\theta_{0,R,\ep} =\Lambda^{s}(w_{0}\chi_{R})* \rho_\ep=\left(\chi_R \Lambda^{s}w_0 \  + \int \frac{w_{0}(y)(\chi_{R}(x)-\chi_{R}(y))}{\vert x-y \vert^{2+s}} \ dy\right) * \rho_\ep.$$
 Since $ \Lambda^{s}w_0 \in L^{\infty}$  then $\chi_R \Lambda^{s}w_0 \in L^{\infty}$. Now, we have to check that the integral in the bracket lies in $L^{\infty}$ uniformly in $R$. For all $R>1$, we split and use the Young inequality $(L^\infty * L^1 \subset L^\infty)$ to get:
\begin{eqnarray*}
&& \int    \frac{w_{0}(y)(\chi_{R}(x)-\chi_{R}(y))}{\vert x-y \vert^{2+s}} \ dy \\ &&=  \int_{\vert x -y \vert < R} \frac{w_{0}(y)(\chi_{R}(x)-\chi_{R}(y))}{\vert x-y \vert^{2+s}} \ dy + \int_{\vert x -y \vert >R} \frac{w_{0}(y)(\chi_{R}(x)-\chi_{R}(y))}{\vert x-y \vert^{2+s}} \ dy \\
 &&\leq  \int_{\vert x-y \vert < R} \frac{ w_{0}(y) \Vert \nabla \chi_R \Vert_{L^\infty}}{\vert x-y \vert^{1+s}} \ dy+2 \int_{\vert x-y \vert >R} \frac{w_{0}(y) \Vert \chi_R \Vert_{L^\infty}}{\vert x-y \vert^{2+s}} \ dy\\
 &&\leq \Vert \nabla \chi_R \Vert_{L^\infty}  \int_{\vert y \vert < R} \frac{ w_{0}(x-y) }{\vert y \vert^{1+s}} \ dy+ 2 \Vert \chi_R \Vert_{L^\infty} \int_{\vert y \vert > R} \frac{w_{0}(x-y) }{\vert y \vert^{2+s}} \ dy.
 \end{eqnarray*}
 We thus get,
 \begin{eqnarray*}
\left\Vert \int \frac{w_{0}(y)(\chi_{R}(x)-\chi_{R}(y))}{\vert x-y \vert^{2+s}} \ dy \right\Vert_{L^{\infty}}
 &\leq& R^{-s} \Vert \nabla \chi \Vert_{L^{\infty} } \Vert w_0 \Vert_{L^{\infty}}   +2 R^{-s} \Vert \chi_R \Vert_{L^{\infty}} \Vert w_0 \Vert_{L^{\infty}} \ 
    \end{eqnarray*}
 Then using the Young inequality and taking the supremum over $R>1$ and $\ep>0$ yields
 $$
\sup_{R>1, \ep>0} \Vert \theta_{0,R,\ep>0} \Vert_{L^{\infty}}= \sup_{R>1, \ep>0}  \Vert  \Lambda^{s}(w_{0}\chi_{R})* \rho_\ep \Vert_{L^{\infty}}  \leq  \sup_{R>1, \ep>0} \Vert \theta_{0,R}  \Vert_{L^{\infty}} \Vert \rho_\ep \Vert_{L^1} \lesssim  \sup_{R>1} \Vert \theta_{0,R}  \Vert_{L^{\infty}} 
 $$
 \qed 
 
 
We also need the following uniform control  on the truncated and regularized  initial data.  

\begin{lemm}
If  $\theta_{0} \in L^\infty$ and $w_{0} \in  H^{s}_{uloc}(\mathbb{R}^2)$ then $\sup_{R>1, \ep>0} \Vert w_{0,R, \ep}  \Vert_{ H^{s}_{uloc}(\mathbb{R}^2)} < \infty$
\end{lemm}
\noindent{\bf Proof.} It is equivalent to show that, for all $R>1$ and $\ep>0$,  $\Lambda^{s}w_{0,R,\ep}  \in L^{2}_{uloc}$.  As in the previous proof, we write
$$
  \Lambda^{s}(w_{0}\chi_{R})*\rho_\ep=\left(\chi_R \Lambda^{s}w_0 \  + C_s \int \frac{w_{0}(y)(\chi_{R}(x)-\chi_{R}(y))}{\vert x-y \vert^{2+s}} \ dy \right)*\rho_{\ep}
$$
Since $w_{0} \in H^{s}_{uloc}(\mathbb{R}^2)$ then  $ \Lambda^{s}w_0 \in L^{2}_{uloc}$ and for all $R>1$,  we see that $\chi_R \Lambda^{s}w_{0,R} \in L^{2}_{uloc}$. \\

\noindent For the second term, we write
\begin{eqnarray*}
  &&\int \frac{w_{0}(y)(\chi_{R}(x)-\chi_{R}(y))}{\vert x-y \vert^{2+s}} \ dy \\ &&=  \int_{\vert x -y \vert < R} \frac{w_{0}(y)(\chi_{R}(x)-\chi_{R}(y))}{\vert x-y \vert^{2+s}} \ dy + \int_{\vert x -y \vert >R} \frac{w_{0}(y)(\chi_{R}(x)-\chi_{R}(y))}{\vert x-y \vert^{2+s}} \ dy \\
 &&\leq   \int_{\vert x-y \vert < R}  \frac{\Vert \nabla \chi_R \Vert_{L^\infty} w_{0}(y)}{\vert x-y \vert^{1+s}} \ dy +  2 \int_{\vert x-y \vert > R} \frac{\Vert \chi_R \Vert_{L^\infty} w_{0}(y)}{\vert x-y \vert^{2+s}} \ dy\\
 &&\leq  R^{-s} \Vert \nabla \chi \Vert_{L^\infty} \int_{\vert y \vert < R}  \frac{ w_{0}(x-y)}{\vert y \vert^{1+s}} \ dy + 2 \Vert \chi_R \Vert_{L^\infty} \int_{\vert y \vert > R} \frac{ w_{0}(x-y)}{\vert y \vert^{2+s}} \ dy. 
 \end{eqnarray*}
 Thus we obtain
 \begin{eqnarray*}
\left \Vert  C_s \int \frac{w_{0}(y)(\chi_{R}(x)-\chi_{R}(y))}{\vert x-y \vert^{2+s}} \ dy \right \Vert_{L^{2}_{uloc}} &\lesssim&  (R^{-s} \Vert \nabla \chi \Vert_{L^\infty} +2 \Vert \chi_R \Vert_{L^\infty} ) \Vert w_0 \Vert_{L^{2}_{uloc}(\mathbb{R}^2)}
 \end{eqnarray*}
where we have use the fact that, since $w_{0} \in H^{s}_{uloc}(\mathbb{R}^2)$ then $w_{0} \in L^{2}_{uloc}(\mathbb{R}^2) $ and $ \mathds{1}_{\vert y \vert < R } \frac{1}{\vert y \vert^{1+s}} \in L^1$  then $\int \frac{w_{0}(y)}{\vert x-y \vert^{1+s}} \ dy \in L^{2}_{uloc}$. \\

\noindent For the second integral, we use that $ \mathds{1}_{\vert y \vert > R}  \frac{1}{\vert y \vert^{2+s}} \in L^1$ and since $w_{0} \in L^{2}_{uloc}$ by convolution we conclude that $ \int_{\vert y \vert > R} \frac{w_{0}(x-y)}{\vert y \vert^{1+s}} \ dy \in L^{2}_{uloc}$. \\ 

\noindent Now, we write 
\begin{eqnarray*}
 \sup_{R>1, \ep>0} \Vert w_{0,R, \ep}  \Vert_{ H^{s}_{uloc}(\mathbb{R}^2)} =   \sup_{R>1, \ep>0} \Vert w_{0,R} * \rho_\ep  \Vert_{ H^{s}_{uloc}(\mathbb{R}^2)} &\leq&  \sup_{R>1, \ep>0}   \Vert w_{0,R}\Vert_{ H^{s}_{uloc}(\mathbb{R}^2)}  \Vert \rho_\ep \Vert_{L^1} \\ &\lesssim& \sup_{R>1}   \Vert w_{0,R}\Vert_{ H^{s}_{uloc}(\mathbb{R}^2)}
 \end{eqnarray*}
 
 \qed


\subsection{Energy inequality and uniform bounds} 

In this section, we prove the following theorem:
\begin{theo}
Suppose that $1/2 < s < 1$, and that $w_{R,\ep}$ is a weak solution of $(SQG)_{R,\ep}$, then  we have the following energy inequality
\begin{equation} \label{eq:equation1}
\partial_{t} A_{\phi}(w_{R,\ep}) \leq C \ \Vert w_{R,\ep} \Vert^{2}_{ H^{s}_{uloc}(\mathbb{R}^2)} 
\end{equation}
where $C$ is a constant depending only on  $\Vert \theta_{0} \Vert_{{L^\infty}(\mathbb{R}^2)}$ and $ \Vert w_{0} \Vert_{ H^{s}_{uloc}(\mathbb{R}^2)}.$
\end{theo}

 \noindent{\bf Proof.}
Let $w_{{R,\ep}}$ be a weak solution to $(SQG)_{{R,\ep}}$ equation.  Recall that, 
$$
\displaystyle \Vert w_{{{R,\ep}}} \Vert^{2}_{{H^{s}_{uloc}} (\mathbb R^2)}  = \sup_{\phi\in B  }\int \left(\frac{\vert w_{{{R,\ep}}} \vert^2}{2}\phi+\frac{\vert \Lambda^{s}w_{{{R,\ep}}}\vert^2}{2}\phi\right) \ dx=\sup_{\phi\in B   } A_{\phi}(w_{{{R,\ep}}}).
$$

\noindent We have
\begin{eqnarray*}
\partial_{t}  A_{\phi}(w_{{R,\ep}}) &=& \int  \partial_{t}  \left(\frac{\vert w_{{R,\ep}} \vert^2}{2}\phi\ +\ \frac{\vert \Lambda^{s}w_{{{R,\ep}}}\vert^2}{2}\phi\right)\ dx \\
&=&\int w_{{R,\ep}}\phi \partial_{t}w_{{{R,\ep}}}\ +\ \phi \Lambda^{s}w_{{{R,\ep}}} \Lambda^{s}(\partial_{t}w_{{{R,\ep}}}) \ dx.
\end{eqnarray*}
We have to justify that the first equality makes sense in $\mathcal{D}'( (0,T] \times \mathbb{R}^2) $. Let us remark that, since $w$ is more regular than $\theta$, it is enough to show that  $\int \phi \theta_{{R,\ep}} \partial_{t} \theta_{R,\ep} \ dx$ makes sense. We begin with the regularity  of $\partial_{t} \theta_{R,\ep}$, using the equation we have
 $$
 \partial_{t} \theta_{R,\ep}=-\nabla. (\theta_{R,\ep} \ u_{R,\ep} ) - \Lambda \theta_{R,\ep}
 $$
 Since $\theta_{R,\ep} \in L^{2}\dot H^{1/2}$ then $\Lambda \theta_{R,\ep} \in L^{2}\dot H^{-1/2}$, moreover, by the continuity of the Riesz transforms on $L^{2}\dot H^{1/2}$  we infer $u_{R,\ep} \in L^{2}\dot H^{1/2}$. We would like to have $ u_R \theta_R \in L^{2} \dot H^{1/2}$.  Since we have
 $$
\Vert u_{R, \ep} \theta_{R, \ep}\Vert_{\dot H^{1/2}} \lesssim \Vert u_{R, \ep} \Vert_{L^\infty} \Vert \theta_{R, \ep} \Vert_{\dot H^{1/2}} +\Vert u_{R, \ep} \Vert_{\dot H^{1/2}} \Vert \theta_{R, \ep} \Vert_{L^\infty},
$$

\noindent it suffices to show that $u_{R,\ep}\in L^\infty L^\infty$.  Indeed, using the Abidi and Hmidi result, we infer that the solutions we obtain are in the space $L^\infty \dot  B^0_{\infty,1}$. Since the Riesz transforms map continuously the  space $L^\infty \dot  B^0_{\infty,1}$ into itself, and thanks to the embedding $L^\infty \dot  B^0_{\infty,1} \hookrightarrow L^\infty L^\infty,$ we get that $u_{R,\ep} \in L^\infty L^\infty$. Therefore, $u_{R, \ep}\Lambda^{s}w_{{{R, \ep}}} \in L^{2}\dot H^{1/2}$, which we had not before regularized since the Riesz transforms are not continuous on the space $L^\infty L^\infty$. Then,  we have $\nabla. (\theta_{R,\ep} \ u_{R,\ep} ) \in L^{2}\dot H^{-1/2}$ and we conclude that $\partial_{t} \theta \in L^{2}\dot H^{-1/2}$. Since $\phi \theta_{R,\ep} \in L^2 \dot H^{1/2},$ then the term $\int \phi \theta_{{R,\ep}} \partial_{t} \theta_{R,\ep} \ dx$ makes sense and we have:
\begin{equation}
 \partial_{t}  A_{\phi}(w_{{R,\ep}})=\int w_{{R,\ep}}\phi \partial_{t}w_{{{R,\ep}}}\ dx +\int \phi \Lambda^{s}w_{{{R,\ep}}} \Lambda^{s}(\partial_{t}w_{{{R,\ep}}}) \ dx.
\end{equation}
 
\noindent Replacing  $\partial_{t}w_{{{R,\ep}}}$ by $-\left(\Lambda^{-s}\nabla \right)\cdot \left(\Lambda^{s}w_{{R,\ep}} (\mathcal{R}^{\perp}{\Lambda^{s} w_{{R,\ep}}})\right) -\Lambda{w_{{R,\ep}}}$ in $(3)$ leads us to the following equality: 
\begin{multline} \label{4t}
\partial_{t} A_{\phi}(w_{{R,\ep}}) =-\int  w_{{R,\ep}}\phi \Lambda^{-s}\nabla.(u\Lambda^{s}w_{{R,\ep}}) \ dx \ - \int w_{{R,\ep}}\phi   \Lambda{w_{{R,\ep}}} \ dx \\ - \int \phi \Lambda^{s}w_{{R,\ep}} \nabla.(u_{{R,\ep}}\Lambda^{s}w_{{R,\ep}})   \ dx 
 -\int \phi \Lambda^{s}w_{{R,\ep}} \Lambda^{s+1}w_{{R,\ep}} \  dx.
\end{multline}
In the estimation of the two first terms, we use the fact that the Riesz transforms $u_{R, \ep}$ are  controlled  in the space $L^{\infty}L^{2}_{uloc}$, more precisely we have the following lemma.
\begin{lemm} \label{bp}
$ u_{R,\ep}$ is  bounded in $L^{2}_t L^{2}_{uloc}$ and we have the following estimate:
$$
\Vert u_{R,\ep} \Vert_{(L^{2}_t L^{2})_{uloc}} \lesssim  \Vert w_{R, \ep} \Vert_{L^{2}_t H^{s}_{uloc}}.
$$
\end{lemm}

\noindent{\bf Proof.} It is equivalent to show that for all $\psi \in B_{\psi}$,  $\psi u_{R, \ep}$ is uniformly  controlled in $L^{2}_t L^2$. We split  $\psi u_{R, \ep}$ as follows:
\begin{eqnarray*}
   \psi u_{R, \ep}=\psi \mathcal{R}^{\perp}\Lambda^{s} w_{R, \ep} &=&\psi \mathcal{R}^{\perp} \Lambda^{s} ( \phi  w_{R, \ep} )+ \psi \mathcal{R}^{\perp}  \Lambda^{s}  (1-\phi) w_{R, \ep} \\
   &=& (1) + (2).
\end{eqnarray*}
Since $w_{R, \ep} \in {(L^{2}H^{s})_{uloc}}$ then $\phi w_{R, \ep} \in L^{2}_t H^{s}$ and $\Lambda^{s} (\phi w_{R, \ep}) \in L^{2}_t L^2$, by the continuity of the Riesz operator on $L^{2}_t L^2$ we deduce that $\mathcal{R}^{\perp} \Lambda^{s}  (\phi  w_{R, \ep})$ is controlled in $L^{2}_t L^2$ and (1) as well. Moreover we have the following inequalities:
\begin{eqnarray*}
\Vert \mathcal{R}^{\perp} \Lambda^{s}  (\phi  w_{R, \ep})  \Vert_{L^{2}_t L^2} \leq \Vert  \Lambda^{s}  (\phi  w_{R, \ep}) \Vert_{L^{2}_t L^2}  &\leq&  \sup_{\phi \in B}  \Vert  \Lambda^{s}  (\phi  w_{R, \ep}) \Vert_{L^{2}_t L^2}\\ &\leq&   \Vert w_{R, \ep} \Vert_{L^{2}_t \dot H^{s}_{uloc}}.
\end{eqnarray*}

Let us denote by $\tilde K$ the kernel of the operator $\mathcal{R}^{\perp}  \Lambda^{s}$.  We have previously seen that this kernel behaves like $\frac{C_s}{\vert x \vert^{2+s}}$ which is in $L^1$ far from the origin, since $w_{R, \ep} \in L^{2}_t {L^{2}_{uloc}}$ we get
 \begin{eqnarray*}
 \Vert  \mathcal{R}^{\perp}  \Lambda^{s}  (1-\phi) w_{R, \ep}  \Vert_{L^{2}_t L^{2}_{uloc}} & \leq& \Vert \mathds{1}_{x\in (supp \phi)^{c}} \tilde K * \phi w_{R, \ep} \Vert_{L^{2}_t L^{2}_{uloc}} \\
&\leq& \Vert \mathds{1}_{x\in (supp \phi)^{c}} \tilde K \Vert_{L^1} \Vert  \phi w_{R, \ep} \Vert_{L^{2}L^{2}_{uloc}} \\
& \lesssim &  \Vert w_{R, \ep} \Vert_{L^{2}H^{s}_{uloc}}.
\end{eqnarray*}
\qed

Let us estimate each of those four terms; to avoid notational burden, we will omit to write the depedence on $\ep$.
For the first one, we would like to write:
\begin{eqnarray*}
 -\int \Lambda^{-s}\nabla.(w_{R}\phi) u_{R}\Lambda^{s}w_{R} \ dx \ dx &\leq& \Vert  \Lambda^{-s}\nabla.(w_{R}\phi) \Vert_{L^{2}_{uloc}} \Vert u_R \Vert_{L^{2}_{uloc}} \Vert \theta_R \Vert_{L^\infty}.  
\end{eqnarray*}
This is done if and only if we have a negative Sobolev regularity on the term $\Lambda^{-s}\nabla.(w_{R}\phi)$ since we cannot control $u_R$ in a space of positive Sobolev regularity (uniformly with respect to the parameters $R$ and $\ep$). Indeed, controlling $u_R$ in a space of positive Sobolev regularity means that we can get a control of the low frequencies of $u_R$ which is hopeless. This is due to the fact that since $\theta_R \in L^\infty L^\infty$ then $u_R \in L^\infty BMO$  and, in particular, we cannot get a control of $\phi u_R$ because $\phi u_R \notin BMO$ ($BMO$ is a space which is defined only modulo constants). Since
 $\phi w_R \in L^{2}\dot H^{s+1/2}$ then $\Lambda^{-s}\nabla.(w_{R}\phi) \in L^{2} \dot H^{2s-1/2} $ we thus need $2s-1/2>0$ that is $1/4 \leq s\leq 1$.  Then, we can write
 \begin{eqnarray*}
 -\int \Lambda^{-s}\nabla.(w_{R}\phi) u_{R}\Lambda^{s}w_{R} \ dx \ dx &\leq& \Vert  \Lambda^{-s}\nabla.(w_{R}\phi) \Vert_{L^{2}_{uloc}} \Vert u_R \Vert_{L^{2}_{uloc}} \Vert \theta_R \Vert_{L^\infty}  \\
 &\lesssim& \Vert \phi w_R \Vert_{\dot H^{1-s}_{uloc}}  \Vert w_{R} \Vert_{H^{s}_{uloc}}  \Vert \theta_{0,R} \Vert_{L^\infty} \\
 &\lesssim& \sup_{\phi \in B_{\phi_0}}  \Vert \phi w_R \Vert_{\dot H^{s}_{uloc}}  \Vert w_{R} \Vert_{ H^{s}_{uloc}}  \Vert \theta_{0,R} \Vert_{L^\infty} \\
 &\lesssim& \Vert  w_R \Vert^{2}_{\dot H^{s}_{uloc}}
\end{eqnarray*}
where we used the embedding $H^{s}\subset H^{1-s}$ for $1/2<s<1$.  \\


\noindent For the second term,  we use the following decay property:
\begin{lemm}
Let $\phi \in \mathcal{D}(\mathbb{R}^2)$ then, for all $x \in \mathbb{R}^2$ we have $\vert \Lambda \phi (x) \vert \leq \frac{C}{\vert x \vert^3}$ where $C$ is a fixed positive constant.
\end{lemm}
\noindent{\bf Proof.}
Let  $B_R$ denote the ball on $\mathbb{R}^2$ of radius $R>0$  centered at the origin. Suppose that $\phi \in \mathcal{D}(\mathbb{R}^2)$ is such that
\begin{equation}
\\\left \{
\aligned
&\phi(x)&=&\ 1\  \mathrm{on}  \ B_R
\\ \nonumber
& \phi(x)&=& \ 0 \ \mathrm{on} \  {B^{c}_{R+1}}
\endaligned
\right.
\end{equation}
If $x \in {B^{c}_{R+2}}$   then  $$\vert \Lambda \phi (x) \vert \leq \int_{y \in B_R} \frac{\vert \phi(y) \vert}{\vert x-y \vert^3} \ dy \leq \frac{C}{\vert x \vert^3}$$ 
If $x \in B_R$, and $y \in {B^{c}_{R+2}}$  it is the same case as before by symmetry. Now, for $x \in B_R$ and $y \in B_R$, we use the following representation of the fractional Laplacian; since $\phi \in \mathcal{S}$ one has
$$
\Lambda \phi = -\frac{1}{2}  \int \frac{\phi(x+y)+\phi(x-y)-2 \phi(x)}{\vert y \vert^{3}} \ dy.
$$
Using a second order Taylor expansion makes us get rid of the singularity at the origin. More precisely, one has
$$
\frac{\phi(x+y)+\phi(x-y)-2 \phi(x)}{\vert y \vert^{3}} \leq \frac{\Vert \nabla^{2}\phi \Vert_{L^\infty}}{\vert y \vert^{}}
$$
Since the right-hand side of the previous inequality is in $L^{1}(B_R)$, we get the desired result. 

\qed

\noindent Using the C\'ordoba and C\'ordoba inequality, an integration by parts and the previous decay property of $\vert \Lambda \phi \vert $,  we obtain:

\begin{eqnarray*}
- \int w_{R}\phi   \Lambda{w_R} \ dx \leq - \int \phi \Lambda(w_{R}^{2}) \ dx &=& - \int \Lambda\phi \ w_{R}^2 \ dx  \\
&\leq& \sum_{k \in \mathbb{Z}^2} \frac{1}{1+\vert k \vert^3} \int_{[-1,1]^2 + k} \vert w_{R}\vert^2 \ dx \\
&\lesssim& \sup_{k \in \mathbb{Z}^2} \int_{[-1,1]^2 + k} \vert w_{R}\vert^2 \ dx \\
& \lesssim&   \Vert w_{R} \Vert^{2}_{H^{s}_{uloc}}.
\end{eqnarray*} 
 
\noindent For the third term we  integrate by parts, and we use the fact that $u$ is divergence free to get

\begin{eqnarray*}
 -\int \phi\Lambda^{s}w_{R} \nabla.(u_{R}\Lambda^{s}w_{R}) \ dx 
 &=& \int \nabla.(\phi\Lambda^{s}w_{R}). u_{{R}}\Lambda^{s}w_{{R}} \ dx\\
 &=&\int  \phi\nabla(\Lambda^{s}w_{R}).u_R \Lambda^{s}w_{R} \ dx\\
 &&+ \int  \Lambda^{s}w_{R} \nabla\phi. u_{R}\Lambda^{s}w_{R}\ dx \\
 &=&  \int \phi \Lambda^{s} w_R \nabla. (u_R \Lambda^{s} w_R) \ dx  + \int \theta_{R}^2 u_{R}.\nabla{\phi} \ dx.
\end{eqnarray*}
Then, lemma \ref{bp} allows us to get
\begin{eqnarray*}
 -\int \phi\Lambda^{s}w_{R} \nabla.(u\Lambda^{s}w_{R}) \ dx =\frac{1}{2}\int \theta_{R}^2 u_{R}.\nabla{\phi} \ dx 
&\lesssim&  \Vert \theta_{R} \Vert_{L^{\infty}(\mathbb{R}^2)} \Vert \nabla \phi \Vert_{L^{\infty}}\Vert  \theta_{R} \Vert_{L^{2}_{uloc}(\mathbb{R}^2)} \Vert u_{R} \Vert_{L^{2}_{uloc}(\mathbb{R}^2)}  \\
&\lesssim& \Vert  w_{R} \Vert^{2}_{ H^{s}_{uloc}(\mathbb{R}^2)}.
\end{eqnarray*}

\begin{rema} {Concerning this third term and more precisely the integration by parts, without regularization, we cannot  integrate by parts since $ \nabla.(\phi\Lambda^{s}w_{R}) \in L^{2}\dot H^{-1/2}$   we must have $u_{{R}}\Lambda^{s}w_{{R}} \in L^{2}\dot H^{1/2}$ which is done only after regularization. }
\end{rema}


\noindent Let us now estimate the last term; we use  Lemma 4 of C\'ordoba and  C\'ordoba applied to $\theta_R$, namely
$$
- \Lambda^{s}w_{R} \ \Lambda^{s+1}w_{R}= -\theta_{R}\Lambda \theta_{R} \leq -\frac{1}{2} \Lambda(\theta_{{R}}^2).
$$
We obtain
\begin{eqnarray*}
-\displaystyle \int \phi \Lambda^{s}w_{R} \ \Lambda^{s+1}w_{R} \ dx \leq  -\int \phi   \Lambda( \theta^{2}_{R}) \ dx&=& -\int \Lambda \phi \   \theta^{2}_{R} \ dx \\
&\leq& \sum_{k \in \mathbb{Z}^2} \frac{1}{1+\vert k \vert^3} \int_{k+[-1,1]^2} \vert \theta_{R} \vert^{2} \ dx \\
&\lesssim&  \Vert \theta_{R} \Vert^{2}_{L^{2}_{uloc}} \\
&\lesssim& \Vert  w_{R} \Vert^{2}_{H^{s}_{uloc}(\mathbb{R}^2)}.
\end{eqnarray*}
Thus, we get the following estimate

$$
-\int \phi \Lambda^{s}w_R \ \Lambda^{s+1}w_R \ dx \leq C \  \Vert w_R \Vert^{2}_{ H^{s}_{uloc}(\mathbb{R}^2)}. 
$$
All these estimates lead to
\begin{equation}
\partial_{t} A_{\phi}(w_{R}) \leq C \ \Vert w_{R} \Vert^{2}_{ H^{s}_{uloc}(\mathbb{R}^2)}.
\end{equation}

\qed
 
\noindent Integrating  $(5)$ in time $s\in [0,t]$ yields
$$
A_{\phi}(w_{R}(x,t))\leq A_{\phi}(w_{0,R}(x))+C \int^{t}_{0} \Vert w_R \Vert^{2}_{ H^{s}_{uloc}(\mathbb{R}^2)} \ ds. 
$$
We now take the supremum over $\phi\in B$ to get
 \begin{eqnarray*}
\Vert w_{R}(x,t) \Vert^{2}_{H^{s}_{uloc}(\mathbb{R}^2)} &\leq& \Vert w_{0,R}(x) \Vert^{2}_{H^{s}_{uloc}(\mathbb{R}^2)} +C  \int^{t}_{0} \Vert w_R(x,s) \Vert^{2}_{ H^{s}_{uloc}(\mathbb{R}^2)} \ ds \\
&\leq& \sup_{R>1} \Vert w_{0,R}(x) \Vert^{2}_{H^{s}_{uloc}(\mathbb{R}^2)} + C \int^{t}_{0} \Vert w_R(x,s) \Vert^{2}_{ H^{s}_{uloc}(\mathbb{R}^2)} \ ds.
\end{eqnarray*}
 By the Gronwall inequality we conclude that
\begin{equation}
\Vert w_{R}(x,t) \Vert^{2}_{L^{\infty}([0,T] ,H^{s}_{uloc}(\mathbb{R}^2))} \leq \sup_{R>1}  \Vert w_{0,R}(x) \Vert^{2}_{H^{s}_{uloc}(\mathbb{R}^2)} e^{Ct}.
\end{equation} 

\noindent Thus $(x,t)\mapsto w_{R}(x,t)$ is globally bounded in the space ${L^{\infty}([0,T] , \dot H^{s}_{uloc}(\mathbb{R}^2))} $. We also need to have a uniform bound for $(x,t)\mapsto w_{R}(x,t)$ in the space $(L^2 \dot H^{s+1/2})_{uloc}$. In other words, we would like to obtain a control of the norm
$$
\Vert w_R \Vert^{2}_{(L^2 \dot H^{s+1/2})_{uloc}}=\Vert \theta_R \Vert^{2}_{(L^2 \dot H^{1/2})_{uloc}} =\sup_{\phi \in B} \int_{0}^{T} \int \phi \vert \Lambda^{1/2} \theta_R \vert^{2} \ dx \ ds.
$$
 We have seen that (see  equality (\ref{4t}))
$$
\partial_{t} A_{\phi}(w_{R}) +\int \phi \theta_R \Lambda \theta_R \ dx \leq C \    \Vert w_{R} \Vert^{2}_{ H^{s}_{uloc}} .
$$
Integrating in time $s\in [0,T]$ leads to
\begin{equation}
A_{\phi}(w_{R}(x,T)) + \int_{0}^{T} \int \phi \theta_R \Lambda \theta_R \ dx \ ds \leq  A_{\phi}(w_{0,R}(x)) + C \int_{0}^{T}  \Vert w_{R} \Vert^{2}_{H^{s}_{uloc}}\ ds.
\end{equation}
The last term in the left-hand side can be rewritten as
\begin{eqnarray*}
\int_{0}^{T} \int \phi \theta_R \Lambda \theta_R \ dx \ ds=\int_{0}^{T} \int \Lambda^{1/2} \theta_R [\Lambda^{1/2}, \phi] \theta_R \ dx \ ds + \int_{0}^{T} \int \phi \vert \Lambda^{1/2} \theta_R \vert^{2} \ dx \ ds,
\end{eqnarray*}
therefore, the inequality $(7)$ becomes
\begin{eqnarray*}
A_{\phi}(w_{R}(x,T)) + \int_{0}^{T} \int \phi \vert \Lambda^{1/2} \theta_R \vert^{2} \ dx \ ds &\lesssim&  A_{\phi}(w_{0,R}(x)) \ +  \int_{0}^{T} C \ \Vert w_{R,\ep} \Vert^{2}_{ H^{s}_{uloc}} \ ds \\ &&+ \left\vert \int_{0}^{T} \int \Lambda^{1/2} \theta_R [\Lambda^{1/2}, \phi] \theta_R \ dx \ ds \right\vert.
\end{eqnarray*}
Then, we split the last integral in the right hand side as follows:
\begin{eqnarray*}
\int_{0}^{T} \int \Lambda^{1/2} \theta_R [\Lambda^{1/2}, \phi] \theta_R \ dx \ ds &=&\int_{0}^{T} \int \psi \Lambda^{1/2} \theta_R [\Lambda^{1/2}, \phi] \theta_R \ dx \ ds \\ &+& \int_{0}^{T} \int (1-\psi)\Lambda^{1/2} \theta_R [\Lambda^{1/2}, \phi] \theta_R \ dx \ ds. \\ &=& (1)+(2).
\end{eqnarray*}
For the first term, by H\"older and Young inequalities we have, for all $\nu >0$ :
\begin{eqnarray*}
\left\vert \int_{0}^{T} \int \psi \Lambda^{1/2} \theta_R [\Lambda^{1/2}, \phi] \theta_R \ dx \ ds \right\vert &\leq& \int_{0}^{T} \Vert \psi \Lambda^{1/2} \theta_R \Vert_{L^2} \Vert \theta_R \Vert_{L^{2}_{uloc}} \ ds \\
& \leq &  \frac{\nu}{2} \int_{0}^{T} \Vert \psi \Lambda^{1/2} \theta_R \Vert^{2}_{L^2} \ ds + \frac{1}{2\nu} \int_{0}^{T} \Vert \theta_R \Vert^{2}_{L^{2}_{uloc}} \ ds. 
\end{eqnarray*} 
For the second term,  we also use H\"older and Young inequalities to obtain that, for all $\eta>0$, we have:
\begin{eqnarray*}
 \int_{0}^{T} \int (1-\psi)\Lambda^{1/2} \theta_R [\Lambda^{1/2}, \phi] \theta_R \ dx \ ds  &=& \int_{0}^{T} \int (1-\psi) \Lambda^{1/2} \theta_R \Lambda^{1/2} (\phi \theta_R) \ dx \ ds \\
& = & \int_{0}^{T} \int \Lambda^{1/2} ((1-\psi) \Lambda^{1/2}\theta_R) \phi \theta_R \ dx \ ds \\
& \leq & \sum_{\vert k \vert > 5} \frac{1}{\vert k \vert^{5/2}} \int_{0}^{T}  \int \phi_k \vert \Lambda^{1/2} \theta_R \vert \vert \phi \theta_R \vert \ dx \ ds \\
& \leq& \frac{\eta}{2}  \sum_{\vert k \vert > 5} \frac{1}{\vert k \vert^{5/2}} \int_{0}^{T}  \int \phi_k \vert \Lambda^{1/2} \theta_R \vert^{2} \ dx \ ds \\ &+&  \frac{1}{2\eta} \int_{0}^{T}    \Vert \theta_R \Vert^{2}_{L^{2}_{uloc}} \ ds,
\end{eqnarray*}
where we have written that $1-\psi=\sum_{\vert k \vert > 5} \phi_k$. Then, we obtain the following inequality
\begin{eqnarray*}
&& A_{\phi}(w_{R}(x,T))+\int_{0}^{T} \int \phi \vert \Lambda^{s+1/2} w_R \vert^{2} \ dx \ ds \lesssim  A_{\phi}(w_{0,R}(x)) \\ &&+ (\frac{1}{2\eta}+ \frac{1}{2\nu})\int_{0}^{T}\Vert \theta_R \Vert^{2}_{L^{2}_{uloc}} \ ds  + (\frac{\eta}{2}+\frac{\nu}{2}) \sup_{\psi \in B_{\psi}} \int_{0}^{T} \int  \psi \vert \Lambda^{1/2} \theta_R \vert^{2} \ dx \ ds.
\end{eqnarray*}
We take the supremum over all $\phi \in B_\phi$ and then  we choose $\nu$ and $\eta$ small enough so that the norm $ \Vert \Lambda^{1/2}\theta_R \Vert^{2}_{(L^2L^{2})_{uloc}}$  appearing in the right hand side is absorbed by that of the left. We obtain
\begin{eqnarray}
 \Vert w_{R} \Vert^{2}_{ H^{s}_{uloc}} + \Vert \Lambda^{1/2}\theta_R \Vert^{2}_{(L^2L^{2})_{uloc}} \lesssim  \Vert w_{0,R}\Vert^{2}_{ H^{s}_{uloc}}+ \int_{0}^{T} \Vert w_{R} \Vert^{2}_{ H^{s}_{uloc}} \ ds.
\end{eqnarray}
By Gronwall's lemma, we recover the following control, for all $T>0$  $$\Vert w_{R} \Vert^{2}_{L^{\infty}([0,T] ,H^{s}_{uloc})} \leq \Vert w_{0,R} \Vert^{2}_{H^{s}_{uloc}} e^{CT}.$$ This readily gives us the desired $(L^2 \dot H^{s+1/2})_{uloc}$ bound. \\

The final step is to show that the  solutions $w_{R,\ep}$ converge weakly to the solutions $w$, and that $w$ are weak solutions of $(SQG)$.

\section{Passage to the limit ($R \rightarrow \infty$, $\ep \rightarrow 0$) }

In this section, we show that the solutions of the truncated $(SQG)_{R,\ep}$ equation tend to the solutions of the classical $(SQG)$ equation, and that the limit has the expected regularity.\\

Let $\eta : (t,x) \mapsto \eta	(t,x) \in \mathcal{D}([0,T] \times  \mathbb{R}^2)$ be a positive test function. Clearly, we have
$$
 \langle\partial_{t}  \theta_{R,\ep},  \eta \rangle = - \langle \theta_{R,\ep}, \partial_{t}  \eta \rangle \underset{R\to+\infty}{\longrightarrow} -  \langle \theta_{\ep}, \partial_{t}  \eta  \rangle.
$$
Thus,
$$
\displaystyle \partial_{t} w_{R,\ep}  \xrightarrow[R\to+\infty, \ \ep \to 0]{} \partial_{t} w \ \ \rm{in} \ \mathcal{D}'([0,T] \times \mathbb{R}^2).
$$
Now, let us focus on the convergence of the non linear term. As usual, weak convergence is not enough to pass to the limit in the non linear term; we need a strong convergence of $\theta_{R,\ep}$ which is given by the  Rellich compactness theorem.

 Since $\theta_{R,\ep}$ is uniformly  bounded with respect to $R$ and $\ep$ in  $(L^{2}\dot H^{1/2})_{loc}$ and $u_{R,\ep}$ is uniformly bounded in the space $(L^2L^2)_{loc} . $ Therefore, by the Sobolev embedding,  the product  $\theta_{R,\ep} u_{R,\ep}$ is uniformly bounded  in the space $(L^2 L^4)_{loc} \cap (L^2L^2)_{loc} \subset (L^2L^{4/3})_{loc} \subset (L^2H^{-1/2})_{loc} $, thus $\nabla. \theta_{R,\ep} u_{R,\ep}$ is uniformly bounded in the space $(L^2 \dot H^{-3/2})_{loc}$.  Moreover, since $\theta_{R,\ep}$ is uniformly bounded in  $(L^{2}\dot H^{1/2})_{loc}$ then $\Lambda \theta_{R,\ep}$ is bounded in $(L^2 \dot H^{-1/2})_{loc}$.  We infer that $\partial_t \theta_{R,\ep}$ is bounded in $(L^2 \dot H^{-3/2})_{loc}$. By Rellich (see e.g \cite{PGLR}, Lemma 21.5, p 218),  we get the strong convergence 
 $$
 \theta_{R,\ep} \xrightarrow[R\to+\infty, \ \ep \to 0]{} \theta \ \text{in} \  (L^{2}L^{2})_{loc}    
 $$
 
  Therefore  $\nabla.(u_{R, \ep} \theta_{R, \ep})$ converges weakly (up to a subsequence) when $R\to+\infty, \ep \to 0$, toward $\nabla.(u \theta)$ in  the space $(L^{2}H^{-3/2})_{loc} $ and thus also in $\mathcal{D}' ((0,T], \mathbb{R}^2)$. 
  Moreover,   $\Lambda \theta_{R, \ep}$ is uniformly bounded in $(L^2 \dot H^{-1/2})_{loc}$,  then $\Lambda \theta_{R,\ep}$ converges weakly (up to a subsequence) when $R\to+\infty, \ep \to 0$ toward $\Lambda w$ in $(L^2 H^{-1/2})_{loc}$ and thus in $\mathcal{D}' ((0,T], \mathbb{R}^2)$. We conclude that, for all $\eta$ in  $\mathcal{D}([0,T] \times  \mathbb{R}^2)$ we have 
  $$
 \langle\eta(t,x),\partial_{t}w \rangle \ = \  \langle\left(\Lambda^{-s}\nabla \right)\cdot \left(\Lambda^{s}w (\mathcal{R}^{\perp}{\Lambda^{s} w})\right) -\Lambda{w}, \eta \rangle
  $$
  The next step is to show that $\theta$ and $w$ have the expected regularity. Let us recall that the spaces $L^{p}_{uloc}$ (with $p>1$) and $H^{s}_{uloc}$ (with $s \in \mathbb{R}$) are respectively the dual spaces of the following separable Banach spaces:
\begin{eqnarray*}
  WL^{p'}(\mathbb{R}^2)&=&\left\{ \theta \in L^{p'}_{loc}(\mathbb{R}^2), \ \sum_{k \in \mathbb{Z}^2} \Vert  \theta \Vert_{L^{p'}(k+[0,1]^2)} <  \infty \right\} \\
  WH^{-s}(\mathbb{R}^2)&=&\left\{ \theta \in H^{-s}_{loc}(\mathbb{R}^2), \ \sum_{k \in \mathbb{Z}^2} \Vert \phi(x-k) \theta \Vert_{H^{-s}(\mathbb{R}^2} < \infty\right\},  
  \end{eqnarray*} 
  where $1/p+1/p'=1.$ \\
  
   By Rellich, we also get the weak convergence (up to a subsequence):
   
   $$
 \theta_{R,\ep} \xrightarrow[R\to+\infty, \ \ep \to 0]{ }  \theta  \ \ \text{weakly}  \ \mathrm{in} \ (L^{2}\dot H^{1/2})_{loc} 
 .$$
   
    Since we have obtained that $w_{R,\ep}$ is uniformly bounded in $(L^\infty  H^s)_{uloc}$ which is the dual space of the Banach space $L^1 W H^{-s}$ we infer that $ w \in \cap_{t<T} (L^{\infty}( [0,t],  H^{s}))_{uloc}$ and then $ \theta \in  \cap_{t<T} (L^{\infty}([0,t], L^{2}))_{uloc}$. \\
   From the uniform bounds,  $w_{R,\ep} \in (L^2 \dot H^{s+1/2})_{uloc}$ which is the dual space of the Banach space $L^2 W \dot H^{-s-1/2}$ therefore we get that $$w \in \cap_{t<T} ((L^{2}( [0,t],\dot H^{s+1/2}))_{uloc}$$ and then $$\theta \in \cap_{t<T} (L^{2}( [0,t], \dot H^{1/2}))_{uloc} $$ \\
   
   For the initial data we have the following convergence:
\begin{eqnarray}
&& \theta_{0,R,\ep} \xrightarrow[R\to+\infty, \ \ep \to 0]{} \theta_{0} \  \ *-\text{weakly}  \ \mathrm{in} \ L^{\infty},  \\
&& w_{0,R,\ep}  \xrightarrow[\hspace{-0,09cm}R\to+\infty, \ \ep \to 0]{} w_{0} \ *-\text{weakly} \ \mathrm{in} \ H^{s}_{uloc}
     \end{eqnarray} 
     Let us show the first convergence. For all $\phi \in L^1$ we have:
     $$
      \langle\theta_{0,R,\ep}, \phi \rangle =  \langle\theta_{0,R}, \phi * \rho_{\ep} \rangle \underset{\ep\to 0}{\longrightarrow}   \langle\theta_{0,R}, \phi  \rangle,
     $$
and thanks to the uniform bounds previously obtained, namely
\begin{eqnarray*}
 \langle\theta_{0,R}, \phi  \rangle&=& \langle\chi_{R} \theta_0, \phi \rangle +  \langle[\Lambda^{s}, w_{0,R}] \chi_R, \phi \rangle \\ &\leq&  \langle\chi_{R} \theta_0, \phi \rangle + R^{-s} \Vert \nabla \chi \Vert_{L^{\infty} } \Vert w_0 \Vert_{L^{\infty}} + 2 R^{-s} \Vert \chi_R \Vert_{L^{\infty}} \Vert w_0 \Vert_{L^{\infty}},
\end{eqnarray*}
we infer that
$$
 \langle\theta_{0,R}, \phi  \rangle \xrightarrow[R\to+\infty]{}  \  \langle\theta_0, \phi \rangle \ \ *-\text{weakly}  \ \mathrm{in} \ L^\infty.
$$

For the convergence $(10)$, we have seen that  $w_{0,R,\ep}$ is uniformly bounded in $ H^{s}_{uloc}$,  therefore,  since $ H^{s}_{uloc}$ is a dual space, we get the following convergence (up to  a subsequence)
$$\hspace{1,5cm} w_{0,R,\ep} \xrightarrow[R\to+\infty, \ \ep \to 0]{}  w_0 \ \ *-\text{weakly}  \ \mathrm{in} \ H^{s}_{uloc}.$$
 and we conclude that 
 $$ \hspace{2,1cm} \theta_{0,R, \ep} \xrightarrow[R\to+\infty, \ \ep \to 0]{}  \theta_{0}\ \ *-\text{weakly}  \ \mathrm{in} \ \Lambda^{s}( H^{s}_{uloc}).$$

\noindent Passing to the weak limit in $(6)$,  we obtain
$$
\Vert w(x,t) \Vert^{2}_{L^\infty  H^{s}_{uloc}(\mathbb{R}^2)} \lesssim  e^{Ct}
$$

\noindent Thus, for $1/2< s < 1$, the solution does not blow up in finite time so we have global existence of weak solutions when $\theta_0 \in \Lambda^{s} ( H^{s}_{uloc}(\mathbb{R}^2)) \cap L^\infty(\mathbb{R}^2)$. \\

\begin{remak} Actually, the energy inequality obtained in $Sect.3.4$ can be improved and we can show that in the range $1/4\leq s \leq 1/2$ we also have global existence this is done in \cite{LZ}.
The case $0<s<1/4$ is little more technical because of the first term among the fourth terms  studied in Sect.3.4, namely
$$
-\int \nabla \Lambda^{-s} (w_R \phi) u_R \theta_R \ dx.
$$
This term is not easy to handle because of the lack of regularity. Nevertheless, in that case, we can show that we have local existence of weak solutions for large initial data in $\Lambda^{s}({ H^{s}_{uloc}(\mathbb{R}^2)})$,  this is also done in the forthcoming work \cite{LZ}. \\
Let us note that the only role of the critical dissipation term $\Lambda \theta$ is to use the result of Abidi and Hmidi \cite{Hm}. This latter ensures us the local boundedness of the Riesz transform in $L^\infty_t L^\infty$ which is a crucial point in our proof.  \\
Finally, let us point out that the uniqueness of those global solutions is open. As it is well known, this kind of proof based on energy and compactness does not provide  uniqueness.

\end{remak}

\vskip0.2cm\noindent{\bf Acknowledgment}:  I thank my advisor Prof. Pierre-Gilles Lemari{\'e}-Rieusset who brought this subject to my attention for his great help and advices.

\end{document}